\documentclass[aos,nameyear,dvips]{arximspdf}
\usepackage{dcolumn}

\doi{10.1214/09-AOS750}
\volume{38}
\issue{3}
\pubyear{2010}
\firstpage{1341}
\lastpage{1362}

\makeatletter

\newcolumntype{d}[1]{D{.}{.}{#1}}

\newtheorem{theorem}{Theorem}

\makeatother

\begin{document}
\begin{frontmatter}

\title{Adjusted empirical likelihood with high-order precision}
\runtitle{Adjusted empirical likelihood}

\begin{aug}
\author[A]{\fnms{Yukun} \snm{Liu}\thanksref{t1}\ead[label=e1]{liuyk1982@yahoo.cn}} and
\author[B]{\fnms{Jiahua} \snm{Chen}\thanksref{t2}\ead[label=e2]{jhchen@stat.ubc.ca}\corref{}}
\runauthor{Y. Liu and J. Chen}
\affiliation{East China Normal University and Nankai University, and
University~of~British~Columbia}
\address[A]{Department of Statistics\\
East China Normal University\\
Shanghai 200241\\
China\\
\printead{e1}}
\address[B]{Department of Statistics\\
University of British Columbia\\
Vancouver, BC, V6T 1Z2\\
Canada\\
\printead{e2}}
\end{aug}

\thankstext{t1}{Supported in part by the NNSF of China Grants 10601026.}

\thankstext{t2}{Supported by the NSERC of Canada and by MITACS.}

\received{\smonth{12} \syear{2008}}
\revised{\smonth{6} \syear{2009}}

%
\begin{abstract}
Empirical likelihood is a popular nonparametric or semi-parametric
statistical method with many nice statistical properties. Yet when the
sample size is small, or the dimension of the accompanying estimating
function is high, the application of the empirical likelihood method
can be hindered by low precision of the chi-square approximation and by
nonexistence of solutions to the estimating equations. In this paper,
we show that the adjusted empirical likelihood is effective at
addressing both problems. With a specific level of adjustment, the
adjusted empirical likelihood achieves the high-order precision of the
Bartlett correction, in addition to the advantage of a guaranteed
solution to the estimating equations. Simulation results indicate that
the confidence regions constructed by the adjusted empirical likelihood
have coverage probabilities comparable to or substantially more
accurate than the original empirical likelihood enhanced by the
Bartlett correction.
\end{abstract}

%

\end{frontmatter}

\section{Introduction}

In applications such as econometrics, statistical finance and
biostatistics, general estimating equations (GEE) in the form $E \{
g(X; \theta)\} = 0$, where $g(x; \theta)$ is a vector-valued function
of the observation vector $x$ and the parameter vector $\theta$, are
often used to define the parameters of interest [Hansen
(\citeyear{H82}), Liang and Zeger (\citeyear{LZ86}), Kitamura and
Stutzer (\citeyear{KS97}) and Imbens, Spady and Johnson
(\citeyear{ISJ98})]. With a semi-parametric setup, scientists run a low
risk of mis-specifying a probability model for the population under
investigation. Particularly when the parameter is over-identified, that
is, when the dimension of $g$ is larger than the dimension of $\theta$,
the generalized moment method (GMM), the empirical likelihood (EL)
method or its variations can be used for statistical inference [Hansen
(\citeyear{H82}), Owen (\citeyear{O88}), Newey and McFadden
(\citeyear{NM94}), Qin and Lawless (\citeyear{QL94}), Imbens
(\citeyear{I97}), Smith (\citeyear{S97}) and Newey and Smith
(\citeyear{NS04})]. Many researchers, however, find that the finite
sample properties of the statistics based on GMM or EL are often very
different from the asymptotic properties at sample sizes common in
applications [Hall and La Scala (\citeyear{HS90}), DiCiccio, Hall and
Romano (\citeyear{DHR91}), Corcoran, Davison and Spady
(\citeyear{CDS95}), Burnside and Eichenbaum (\citeyear{BE96}), Corcoran
(\citeyear{C98}) and Tsao (\citeyear{T04})]. High-order approximations
to the finite sample distribution based on the Bartlett correction or
bootstrapping can be helpful [DiCiccio, Hall and Romano
(\citeyear{DHR91}), Hall and Horowitz (\citeyear{HH96}), Brown and
Newey (\citeyear{BN02}), Newey and Smith (\citeyear{NS04}) and Chen and
Cui (\citeyear{CC07})]. Yet they do not always live up to their
promise, particularly for high-dimensional data [Corcoran, Davison and
Spady (\citeyear{CDS95}) and Tsao (\citeyear{T04})].

We propose a novel approach via adjusted empirical likelihood (AEL)
[Chen, Variyath and Abraham (\citeyear{CVA08})] to achieve the
high-order precision promised by the Bartlett correction. The AEL is
obtained by adding a pseudo-observation into the data set. Its
principal utility is to overcome the difficulty arising when the
estimating equations have no solution; a solution is required in the EL
approach. By using a conventional level of adjustment, Chen, Variyath
and Abraham (\citeyear{CVA08}) found the AEL improves the approximation
precision of the chi-square limiting distribution. More recently,
Emerson and Owen (\citeyear{EO09}) discussed the level of adjustment
for inference on multivariate population mean. However, the optimal
level of adjustment remains unknown. In this paper, we derive a
high-order expansion of the adjusted empirical likelihood ratio
statistic, specify an optimal level of adjustment that enables the
high-order approximation, prove that the resulting AEL shares the same
high-order precision as the Bartlett corrected EL (BEL) and construct a
less biased estimator of the Bartlett correction factor that
effectively improves the approximation precision.

Although the AEL and the BEL have the same high-order precision, their
finite sample performances differ. Simulation studies show that the AEL
has better precision than the BEL in general, and especially under
linear and asset-pricing models. The AEL with conventional level of
adjustment, AEL$_0$, is found to have comparable precisions to the AEL
under many models considered, but it lacks some generality. In
particular, the AEL improves over the AEL$_0$ under linear and
asset-pricing models.

\section{The EL and the Bartlett correction}
\subsection{The empirical likelihood}
To convey the idea, suppose we have $x_1, x_2,\break \ldots, x_n$ as a
random sample from a nonparametric population $F(x)$ such that $x
\in\mathbb{R}^m$ with dimension $m$. Assume that the GEE model is
defined by
\[
E g(X; \theta) = 0
\]
for a $q$-dimensional estimating function $g$ and a $p$-dimensional
parameter $\theta$. The profile empirical likelihood function of
$\theta$ is defined as
%
%
\begin{equation}\label{EL.def}
L_n(\theta) = \sup \Biggl\{ \prod_{i=1}^n p_i\dvtx p_i \geq0, \sum_{i=1}^n
p_i=1, \sum_{i=1}^n p_i g(x_i; \theta)=0 \Biggr\}.
\end{equation}
The empirical log-likelihood ratio function is defined by $R_n(\theta)
= -2 \log(n^n\times\break L_n(\theta))$ [see Owen (\citeyear{O01}) and Qin and
Lawless (\citeyear{QL94})]. One celebrated property of the empirical
likelihood is that under some general conditions,
\[
\mbox{\textsc{pr}}\{ R_n(\theta_0) \leq x \} = \mbox{\textsc{pr}}\{
\chi_q^2 \leq x \} + O(n^{-1})
\]
as $n \to\infty$ where $\theta_0$ is the true parameter value. This
property is most convenient for the construction of confidence regions
of $\theta$,
%
%
\begin{equation}\label{CI}
\{ \theta\dvtx R_n(\theta) \leq c(1-\alpha; q) \}
\end{equation}
with $c(1-\alpha; q)$ being the $(1-\alpha)$th quantile of the
chi-square distribution with $q$ degrees of freedom, and $1-\alpha$
being the pre-selected confidence level. Such confidence regions are
renowned for their data-driven shape, and there is no need to estimate
any scalar parameters. For other results, such as when $\theta_0$ is
replaced by its nonparametric maximum EL estimate $\hat\theta$, we
refer to Qin and Lawless (\citeyear{QL94}).

\subsection{The Bartlett correction of the EL}\label{sec22}
The precision of the confidence region constructed by (\ref{CI}) can be
poor, particularly when the sample size is small. To improve the
precision of the coverage probability, we may calibrate the
distribution of $R_n(\theta_0)$ by bootstrapping or by high-order
approximations. We now review high-order approximation via the Bartlett
correction.

The Bartlett correction for a smooth function of means was first
established by DiCiccio, Hall and Romano (\citeyear{DHR91}) while
estimating questions by Chen and Cui (\citeyear{CC06},
\citeyear{CC07}). For ease of illustration, we consider the situation
where $p = q=1$ and $g(x; \theta) = x - \theta$. Under this model, the
parameter $\theta$ is the population mean. The chi-square approximation
has precision $O(n^{-1})$ and the confidence interval of $\theta$ based
on the chi-square approximation may not have accurate coverage
probabilities. The Bartlett correction can improve the approximation
precision to $O(n^{-2})$.

By the Lagrange method, when the solution to ${\sum_{i=1}^{n}}p_i
g(x_i; \theta) = 0$ exists, we have
\[
R_n(\theta) = {\sum_{i=1}^{n}}\log\{ 1 + \lambda g(x_i; \theta) \}
\]
for a Lagrange multiplier $\lambda$ that is the solution to
%
%
\begin{equation}
\label{eqn4}
\sum_{i=1}^{n}\frac{ g(x_i; \theta) }{ 1 + \lambda g(x_i;
\theta)} = 0.
\end{equation}
Let $\alpha_r = E\{ g (X; \theta)\}^r$ and $A_r = n^{-1}
{\sum_{i=1}^{n}} \{g(x_i; \theta)\}^r - \alpha_r$. Without loss of
generality, we assume that either $\alpha_2=1$ or we can replace $g(x;
\theta)$ with $\alpha_2^{-1/2} g(x; \theta)$. Assuming that $\theta$ is
the true parameter value, we can write
\[
\lambda= \lambda_1 + \lambda_2 + \lambda_3 + O_p(n^{-2})
\]
with
\begin{eqnarray*}
\lambda_1 &=& A_1,\qquad \lambda_2 = \alpha_3 A_1^2 - A_1 A_2,
\\
\lambda_3 &=& A_1 A_2^2 + A_1^2 A_3 + 2 \alpha_3^2 A_1^3 - 3
\alpha_3A_1^2 A_2 - \alpha_4 A_1^3.
\end{eqnarray*}
Under some moment conditions, $\lambda_r = O_p(n^{-r/2})$ for $r = 1,
2, 3$. Substituting these expansions into the expression for
$R_n(\theta)$, we get
%
%
\begin{equation}
R_n(\theta) = n \{ R_1 + R_2 + R_3\}^2 + O_p(n^{-3/2})
\end{equation}
with
\begin{eqnarray*}
R_1 &=&  A_1, \\
R_2 &=& \tfrac{1}{3} \alpha_3 A_1^2 - \tfrac{1}{2} A_1 A_2,
\\
R_3 &=& \tfrac{3}{8} A_1 A_2^2 + \tfrac{4}{9} \alpha_3^2 A_1^3
- \tfrac{5}{6} \alpha_3 A_1^2 A_2 + \tfrac{1}{3} A_1^2 A_3
- \tfrac{1}{4} \alpha_4 A_1^3.
\end{eqnarray*}
DiCiccio, Hall and Romano (\citeyear{DHR91}) find that
the cumulants of
\[
n(1- b/n)(R_1 + R_2 + R_3)^2
\]
match those of the $\chi_1^2$ distribution to the order of
$n^{-3/2}$ when
%
%
\begin{equation}\label{Bartlett}
b=\tfrac{1}{2} \alpha_4 - \tfrac{1}{3} \alpha_3^2.
\end{equation}
Furthermore, since $R_1+R_2+R_3$ are smooth functions of general sample
means, the result of Bhattacharya and Ghosh (\citeyear{BG78}) implies
that
\[
\mbox{\textsc{pr}}\{ n(1- b/n)(R_1 + R_2 + R_3)^2 \leq x \}
= \mbox{\textsc{pr}}\{ \chi_1^2 \leq x \} + O(n^{-2}).
\]
More details are in the \hyperref[app]{Appendix}.

In applications, the value $b$ must be replaced by some root-$n$
consistent estimator, and, in theory, the replacement does not affect
the high-order asymptotic conclusion. Naturally, $b$ is often replaced
by a moment estimate.

Another way to improve the finite sample performance is to use
bootstrap calibration, that is, to estimate the sample distribution of
the $R_n(\theta)$ via a bootstrap resampling scheme [see, e.g., Hall
and Horowitz (\citeyear{HH96})]. There are situations where the
solution $p_i$'s to the constraints in (\ref{EL.def}) at $\theta=
\theta_0$ do not exist with nonnegligible probability. A convention
adopted in this situation is to define $R_n(\theta) = \infty$. However,
if ${\mbox{\textsc{pr}}}\{R_n(\theta_0) = \infty\} > \alpha$, then
${\mbox{\textsc{pr}}}\{ R_n(\theta_0) < c \} < 1 - \alpha$ for any
finite $c$. Consequently, a bootstrap scheme can at most boost the
coverage probability to $1 - {\mbox{\textsc{pr}}}\{R_n(\theta_0) =
\infty\}$ which is still below the nominal level $1 - \alpha$. This
problem is clearly also shared by the Bartlett correction [see also
Tsao (\citeyear{T04})].

\section{The AEL and the high-order approximation}
\subsection{The adjusted empirical likelihood}
For each given $\theta$, the likelihood ratio function
$R_n(\theta)$ is well defined only if the convex hull of
%
%
\begin{equation}\label{hull}
\{ g(x_i; \theta)\dvtx i=1, 2, \ldots, n \}
\end{equation}
contains the $q$-dimensional vector $\mathbf{0}$. When $n$ is not
large, or when a good candidate (vector) value of $\theta$ is not
available, this convex hull often fails to contain $\mathbf{0}$ [see,
e.g., Chen, Variyath and Abraham (\citeyear {CVA08})]. Blindly setting
$L_n(\theta) = 0$ as suggested in the literature fails to provide
information on whether $\theta$ is grossly unfit to the data or is in
fact only slightly off an appropriate value. Let $g_i = g(x_i;
\theta)$, $i=1, \ldots, n$, and
\[
g_{n+1} = - a_n \bar g_n = - a_n n^{-1} {\sum_{i=1}^{n}}g_i
\]
for some $a_n > 0$. The adjusted (profile) empirical likelihood is
defined as
%
%
\begin{equation}\label{AEL.def}
L_n(\theta; a_n) = \sup \Biggl\{ \prod_{i=1}^{n+1} p_i\dvtx p_i \geq0,
\sum_{i=1}^{n+1} p_i = 1, \sum_{i=1}^{n+1} p_i g_i = 0 \Biggr\}
\end{equation}
and the adjusted empirical likelihood ratio function as
\[
R_n(\theta; a_n) = - 2 \log\{ (n+1)^{n+1} L_n(\theta; a_n) \}.
\]
Because $\bar g_n$ and $g_{n+1}$ are on opposite sides of $\mathbf
{0}$, the AEL is always well defined. Namely, its value is always
nonzero. When $a_n = o_p(n^{2/3})$, Chen, Variyath and Abraham
(\citeyear{CVA08}) showed that the first-order asymptotic properties of
the EL are retained by the AEL, and a conventional $a_n=\max\{1, \log
n/2\}$ was found useful in a number of examples. However, an optimal
choice of $a_n$ remains unsolved. We next recommend a specific $a_n$
and show that the resulting AEL achieves the goal attained by the
Bartlett correction.

\subsection{AEL with high-order precision}

The level of adjustment at which the AEL has high-order precision is
$a_n=b/2$, where $b$ is the Bartlett correction factor for the usual
EL. This surprising relationship reveals an intrinsic relationship
between the AEL and the Bartlett correction. Indeed, the proof of the
following result is built on the Bartlett correction.
\begin{theorem}
\label{thm1} Suppose that $x_1, x_2, \ldots, x_n$ is a random sample
from an $m$-variate nonparametric population $F(x)$. Assume that the
GEE model is defined by
\[
E g(X; \theta) = 0,
\]
where $\theta$ is a $p$-dimensional parameter, $g(X; \theta)$ is a
$q$-dimensional estimating function, and its characteristic function
satisfies Cram\'{e}r's condition,
\[
{\limsup_{\| t \| \rightarrow\infty}}
| E \exp\{\mathbf{i}t^T g(X; \theta) \} |<1.
\]
Assume also that $E\Vert g(X; \theta) \Vert^{18} < \infty$
and $\operatorname{var}( g(X; \theta))$ is positive definite.

Let $\theta_0$ be the true parameter value and $a_n = a +
O_p(n^{-1/2})$. Then
\[
R_n(\theta_0; a_n) = n \{R_1+R_2 + R_{3a}\}^T \{R_1+R_2 + R_{3a}\} +
O_p(n^{-3/2}),
\]
where $R_1$, $R_2$ and $R_{3a}$ will be given in (\ref{R1}) and
(\ref{Ra}). When $a = b/2$ where $b$ is the Bartlett correction factor
for the usual empirical likelihood,
\[
\mbox{\textsc{pr}} \bigl\{ n \{R_1+R_2 + R_{3a}\}^T \{R_1+R_2 + R_{3a}\}
\leq x \bigr\} = \mbox{\textsc{pr}}( \chi_q^2 \leq x) + O(n^{-2}).
\]
\end{theorem}

Adding a pseudo-observation $g_{n+1}$ results in a slightly different
$R_{3a}$ as compared to $R_3$ in Section \ref{sec22}. This explains the
choice of the notation.

When $q=1$, the Bartlett correction factor $b = \alpha_4/2 -
\alpha_3^2/3 > 0$ unless $g(X; \theta)$ degenerates. Hence, the
pseudo-observation obtained by setting $a_n = b/2$ or its suitable
estimator satisfies the condition $a_n > 0$ required by the AEL. When
$q > 1$, it is uncertain whether $b > 0$ or not. While Theorem
\ref{thm1} remains valid, there is a small probability that the AEL is
not defined when $b < 0$. We can easily avoid this problem by adding
two pseudo-observations. Let
%
%
\begin{equation}\label{AEL.def2}
L_n(\theta; a_{1n}, a_{2n}) = \sup \Biggl\{\prod_{i=1}^{n+2} p_i\dvtx p_i
\geq0, \sum_{i=1}^{n+2} p_i = 1, \sum_{i=1}^{n+2} p_i g_i = 0 \Biggr\}
\end{equation}
and let the adjusted empirical likelihood ratio function be
\[
R_n(\theta; a_{1n}, a_{2n})
= - 2 \log\{ (n+2)^{n+2} L_n(\theta; a_{1n}, a_{2n}) \}
\]
with $g_{n+1} = - a_{1n} \bar g$ and $g_{n+2} = a_{2n} \bar g$. When
$a_{2n} - a_{1n} = b$, the result of Theorem~\ref{thm1} remains.

In general, the Bartlett correction factor $b$ can be written as the
difference of two positive values. This decomposition gives us natural
choices of $a_{1n}$ and $a_{2n}$ for multidimensional estimating
functions. In simulations, we added a single pseudo-observation when
$q=1$ and two pseudo-observations when $q\geq2$. We also recommend this
practice in applications. More detailed discussions about the Bartlett
correction factor $b$ are given in the next subsection.

When $q > p$ where the parameter is over-identified, it is more
efficient to construct confidence regions with
\[
\Delta_n(\theta; a_n) = R_n(\theta; a_n) - \inf_\theta R_n(\theta; a_n).
\]
When $a_n = 0$, Chen and Cui (\citeyear{CC07}) show that
$\Delta_n(\theta_0; 0)$ is also Bartlett correctable. The result of
Theorem \ref{thm1} remains valid as follows.
\begin{theorem}
\label{thm2} Assume the same conditions as in Theorem \ref{thm1}, and
that there exists a neighborhood of $\theta_0$, $N(\theta_0)$ and an
integrable function, $h(x)$, such that
\[
{\sup_{\theta\in N(\theta_0)}}\| \partial^3 g(x; \theta)/\partial
\theta^3\|^3 \leq h(x).
\]
Then at the level of adjustment $a_n = a + O_p(n^{-1/2})$,
\[
\Delta_n(\theta_0; a_n) = n \{R_1+R_2 + R_{3a}\}^T \{R_1+R_2 + R_{3a}\}
+ O_p(n^{-3/2})
\]
for some $R_1$, $R_2$ and $R_{3a}$,
and there exists a Bartlett correction factor $b$ such that when
$a = b/2$,
\[
\mbox{\textsc{pr}} \bigl\{ n \{R_1+R_2 + R_{3a}\}^T \{R_1+R_2 +
R_{3a}\} \leq x \bigr\} = \mbox{\textsc{pr}}({ \chi_p^2} \leq x) +
O(n^{-2}).
\]
\end{theorem}

The expressions of the Bartlett correction factor $b$ and $R_j$, $j=1,
2$, in Theorem~\ref{thm2} are the same as in Chen and Cui
(\citeyear{CC07}). When $a_n = 0$, $R_{3a}$ also becomes their $R_3$.
More details and a brief proof are given in the
\hyperref[app]{Appendix}.

\subsection{Estimation of the Bartlett correction factor $b$}
We first consider the estimation of $b$ in the case of Theorem
\ref{thm1}. Even for the simplistic one-sample problem,
Bartlett-corrected ordinary EL confidence intervals for the population
mean often have lower than nominal coverage probabilities when the
Bartlett correction factor $b$ is replaced by its moment estimator. The
Bartlett-corrected EL intervals with theoretical $b$ are often much
more satisfactory. Our investigation reveals that the moment estimator
of $b$ usually grossly under estimates particularly when $n$ is small,
say $n=20$, 30. See the simulation results presented in the next
section.

Let us first examine the case of $q = p = 1$ where the Bartlett
correction factor is given by
\[
b=\frac{\alpha_4}{2\alpha_2^2} - \frac{\alpha_3^2}{3\alpha_2^3}.
\]
Note that we no longer assume $\alpha_2 = 1$. The moment estimators of
$\alpha_r$ are given by $\hat\alpha_r = n^{-1} \sum_{i=1}^n (g_i - \bar
g)^r$. Since $E \hat\alpha_2 = (n-1)\alpha_2/n $, we estimate
$\alpha_2$ by $\tilde\alpha_2 = n \hat\alpha_2 /(n-1)$ to reduce bias.
In summary, we use the estimators given in the following table to
construct a less-biased estimator of $b$:
\[
\matrix{
\mbox{Parameter} & \mbox{Estimator} & \mbox{Expression}\cr
\alpha_2 & \tilde\alpha_2 & n\hat\alpha_2/(n-1)\cr
\alpha_4 & \tilde\alpha_4 & (n \hat{\alpha}_4 -6 \tilde
{\alpha}_2^2)/ (n-4)\cr
\alpha_2^2 & \tilde\alpha_{22} &
\tilde{\alpha}_2^2 - \tilde{\alpha}_4/n\cr
\alpha_3 & \tilde\alpha_3 & n\hat\alpha_3/(n-3) \cr
\alpha_3^2 & \tilde\alpha_{33} &
\tilde{\alpha}_3^2 - (\hat{\alpha}_6-\tilde\alpha_3^2)/n\cr
\alpha_2^3 & \tilde\alpha_{222} & \tilde\alpha_2^3 .}
\]
The above choices are motivated as follows. Since
\begin{eqnarray*}
E \hat{\alpha}_4 &=& \alpha_4 - \frac{4 \alpha_4}{n} +
\frac{6\alpha_2^2}{n} +
O(n^{-2}), \\
E \tilde\alpha_2^2
&=& \alpha_2^2 + \frac{\alpha_4}{n} + O(n^{-2}),\\
E \hat\alpha_3 &=& \alpha_3 - \frac{3\alpha_3}{n} +O(n^{-2}),
\end{eqnarray*}
we estimate $\alpha_4$, $\alpha_2^2$ and $\alpha_3$ by $ \tilde\alpha_4
= (n \hat{\alpha}_4 -6 \tilde{\alpha}_2^2)/(n-4), $ $ \tilde\alpha_{22}
= \tilde{\alpha}_2^2 -\tilde{\alpha}_4/n $ and $\tilde\alpha_3 =
n\hat\alpha_3/(n-3) $, respectively. The biases of $\tilde\alpha_4$,
$\tilde\alpha_{22}$ and $\tilde\alpha_{3}$ are of order $O(n^{-2})$
compared to the $O(n^{-1})$ biases of the corresponding moment
estimators. Precise form of the $O(n^{-1})$ bias of $\hat\alpha^2_3$ is
complex. Hence, we aim to reduce rather than completely eliminate the
$O(n^{-1})$ bias. Since $\tilde\alpha_3 \approx\frac{1}{n}\sum_{i=1}^n
g_i^3$, we have approximately $E\tilde\alpha_3 = \alpha_3 $ and $ E
\tilde\alpha_3^2 = \alpha_3 ^2 + \operatorname{var}(\tilde\alpha_3)$,
and approximately $ \operatorname{var}(\tilde\alpha_3) =
(\alpha_6-\alpha_3^2)/n. $

When $q = p > 1$, the expression for $b$ is more complex. Let
$V(\theta) =\break \operatorname{var}\{g(X; \theta)\}$ be the covariance
matrix. By eigenvalue decomposition, we may write
\[
V(\theta_0) = P \operatorname{diag}\{ \xi_1, \ldots, \xi_q \} P^T
\]
such that $PP^T = I$ and $\xi_1, \ldots, \xi_q$ are eigenvalues of
$V(\theta_0)$. Furthermore, let $Y=P^T g(X; \theta_0)$, and for any
positive integers $(r, s, \ldots, t)$, define
%
%
\begin{eqnarray}
\alpha^{r s \cdots t} = E\{Y^{r}Y^{s}\cdots Y^{t}\},
\end{eqnarray}
where $Y^t$ is the $t$th component of vector $Y$.

It can be seen that after $g$ is transformed by multiplying $P$,
$\alpha^{rr} = \xi_r$ and $\alpha^{rs} = 0$ for $r \neq s$. The
Bartlett correction factor can then be written as
\begin{eqnarray*}
b &=& \frac{1}{q} \biggl\{ \frac{1}{2} \sum_{r, s}
\frac{\alpha^{rrss}}{\alpha^{rr} \alpha^{ss}} - \frac{1}{3} \sum_{r, s,
t} \frac{\alpha^{rst} \alpha^{rst}}
{ \alpha^{rr} \alpha^{ss} \alpha^{tt}} \biggr\}\\
&=& \frac{1}{q} \biggl\{ \sum_{r} \frac{\alpha^{rrrr}}{2(\alpha^{rr} )^2} +
\sum_{r \neq s}\frac{\alpha^{rrss}}{ 2\alpha^{rr} \alpha^{ss}}
\biggr\} \\
&&{} - \frac{1}{q} \biggl\{ \sum_{r} \frac{(\alpha^{rrr})^2}{3(\alpha^{rr} )^3}
+ \sum_{r \neq s} \frac{(\alpha^{rss})^2}{\alpha^{rr}(\alpha^{ss} )^2}
+ 2\sum_{r < s < t}
\frac{(\alpha^{rst})^2}{\alpha^{rr}\alpha^{ss}\alpha^{tt}} \biggr\}
\\
&=& \frac{1}{q}\sum_{r} \biggl\{ \frac{\alpha^{rrrr}}{2(\alpha^{rr} )^2} -
\frac{(\alpha^{rrr})^2}{3(\alpha^{rr} )^3} \biggr\} + \frac{1}{2q} \sum_{r
\neq s} \biggl\{ \frac{\alpha^{rrss}}{ \alpha^{rr} \alpha^{ss}} -
\frac{(\alpha^{rss})^2}{\alpha^{rr}(\alpha^{ss} )^2}
\biggr\} \\
&&{} - \frac{1}{q} \biggl\{ \frac{1}{2} \sum_{r \neq s}
\frac{(\alpha^{rss})^2}{\alpha^{rr}(\alpha^{ss} )^2} + 2\sum_{r < s <
t} \frac{(\alpha^{rst})^2}{\alpha^{rr}\alpha^{ss}\alpha^{tt}} \biggr\}.
\end{eqnarray*}
Let
\begin{eqnarray*}
b_1 &= & \frac{1}{q}\sum_{r} \biggl\{ \frac{\alpha^{rrrr}}{2(\alpha^{rr} )^2}
- \frac{(\alpha^{rrr})^2}{3(\alpha^{rr} )^3} \biggr\} + \frac{1}{q}\sum_{r <
s} \biggl\{ \frac{\alpha^{rrss}}{ \alpha^{rr} \alpha^{ss}} -
\frac{(\alpha^{rss})^2}{\alpha^{rr}(\alpha^{ss} )^2}
\biggr\}, \\
b_2 &= & \frac{1}{q} \sum_{r < s}
\frac{(\alpha^{rss})^2}{\alpha^{rr}(\alpha^{ss} )^2} +
\frac{2}{q}\sum_{r < s < t}
\frac{(\alpha^{rst})^2}{\alpha^{rr}\alpha^{ss}\alpha^{tt}}.
\end{eqnarray*}
Clearly, both $b_1$ and $b_2$ are positive and $b = b_1 - b_2$. There
can be other ways to decompose $b$. We have chosen the above
decomposition so that both $b_1$ and $b_2$ are of moderate size.

Note that the Bartlett correction factor(s) depends on the unknown
$\theta_0$. In applications, we first compute a maximum adjusted
empirical likelihood estimate $\hat\theta$ at $a_n = \log n/2$, and use
it as a tentative replacement of $\theta_0$ for estimating $b$ or $b_1$
and~$b_2$. We decompose the sample variance of $g(x; \theta)$ at
$\theta= \hat\theta$ to obtain the orthogonal matrix $P$. We then
obtain $Y_i = P^T g(X_i; \hat\theta)$ and define the moment estimators
as
%
%
\begin{eqnarray}
\label{alpha}
\hat\alpha^{rs\cdots t}
=n^{-1} \sum_{i=1}^n Y_i^{r}Y_i^{s}\cdots Y_i^{t}.
\end{eqnarray}

To reduce the bias in the estimation of $b_1$ and $b_2$,
we use the estimators given in the following table:
\[
\begin{array}{c@{\qquad}c@{\qquad}c}
\mbox{Parameter} & \mbox{Estimator} & \mbox{Expression}\\
\alpha^{rr} & \tilde\alpha^{rr} & n\hat\alpha^{rr}/(n-1)\\
\alpha^{rrss} & \tilde\alpha^{rrss} &
\{n \hat\alpha^{rrss} -2 \tilde\alpha^{rr}\tilde\alpha^{ss}
- 4I(r=s){\tilde\alpha^{rr}\tilde\alpha^{rr}}\}/(n-4)
\\
 \alpha^{rst}&\tilde\alpha^{rst} &
n\hat\alpha^{rst}/(n-3) \cr
\alpha^{rst}\alpha^{rst}&\tilde\alpha^{rst,rst}&
\tilde\alpha^{rst}\tilde\alpha^{rst} - (\hat\alpha^{rrsstt}
- \tilde\alpha^{rst}\tilde\alpha^{rst} )/n
\\
\alpha^{rr}\alpha^{ss}&\tilde\alpha^{rr,ss}&
{ \tilde\alpha^{rr}\tilde\alpha^{ss}} - \tilde\alpha^{rrss}/n\cr
\alpha^{rr}\alpha^{ss}\alpha^{tt}&\tilde\alpha^{rr,ss,tt}&
\tilde\alpha^{rr}\tilde\alpha^{ss}\tilde\alpha^{tt}
\end{array}
\]
for all $1\leq r, s,t\leq q$, and $I(r=s)$ is the indicator function.
We denote the resulting estimates as $\tilde b_1$ and $\tilde b_2$. For
$q > 1$, we add two pseudo-observations with $a_{1n} = \tilde b_1/2$
and $a_{2n} = \tilde b_2/2$ in the simulations.

To examine the bias properties of the new estimator, we generated
10,000 sets of random samples from a number of selected univariate,
bivariate and trivariate distributions. The population distributions
are not important at this stage, and they will be specified in the
simulation section. We computed the Bartlett correction factors and
their average estimates for constructing confidence regions of the
population mean. The outcomes are given in Tables \ref{table2} and
\ref{multi2}. The moment estimators are denoted as $ b_n$ and the new
estimators as $\tilde b_n$. Clearly, the new estimators are much less
biased under the normal, exponential and chi-square distributions.
Under mixture distributions, $\tilde b_n$ overestimates $b$, but the
resulting AEL confidence intervals still have good coverage properties.
We also examined the bias properties under a number of linear models.
The results are given in Table \ref{table4}. Again, $\tilde b_n$ is
much less biased. The model specifications are relegated to the
simulation section.

%
%
\begin{table}
\caption{Bartlett correction factors and their average estimates for
univariate population mean}
\label{table2}
\begin{tabular*}{\tablewidth}{@{\extracolsep{\fill}}lccccc@{}}
\hline
$\bolds{n}$ & & $\bolds{N(0,1)}$ & \textbf{Exp(1)} & $\bolds{0.2 N_1 + 0.8 N_2}$
& $\bolds{\chi^2_1}$ \\
\hline
& $b$ & 1.50 & 3.17 & 1.11 & 4.83 \\
$20$ & $b_n$ & 1.16 & 1.40 & 1.14 & 1.59 \\
& $\tilde b_n$ & 1.57 & 3.19 & 2.08 & 5.56 \\
$30$ & $b_n$ & 1.26 & 1.66 & 1.15 & 1.96 \\
& $\tilde b_n$ & 1.56 & 3.17 & 1.63 & 5.12 \\
\hline
\end{tabular*}
\end{table}

%
%
\begin{table}[b]
\caption{Bartlett correction factors and their average estimates for
multivariate ($q=2, 3$) population mean}\label{multi2}
\begin{tabular*}{\tablewidth}{@{\extracolsep{\fill}}lccccc@{}}
\hline
$\bolds{n}$ & & \textbf{(a)} & \textbf{(b)} & \textbf{(c)} & \textbf{(d)} \\
\hline
$q=2$ & $b$ & 3.21 & 3.71 & 1.68 & 2.21 \\
\quad$20$ & $b_n$ & 1.63 & 1.67 & 1.48 & 1.46 \\
& $\tilde b_n$ & 2.93 & 3.34 & 2.55 & 2.14 \\
\quad$30$ & $b_n$ & 1.90 & 1.98 & 1.56 & 1.64 \\
& $\tilde b_n$ & 3.06 & 3.47 & 2.18 & 2.20 \\
[4pt]
$q=3$ & $b$ & 4.07 & 3.84 & 2.36 & 2.67 \\
\quad$30$ & $b_n$ & 2.27 & 2.24 & 1.98 & 2.00 \\
& $\tilde b_n$ & 3.72 & 3.47 & 2.62 & 2.62 \\
\quad$50$ & $b_n$ & 2.67 & 2.61 & 2.13 & 2.22 \\
& $\tilde b_n$ & 3.89 & 3.64 & 2.52 & 2.67 \\
\hline
\end{tabular*}
\end{table}


When $q > p$, we prefer $\Delta_n(a)$ for constructing confidence
intervals as in Theorem \ref{thm2}. However, as indicated in Chen and
Cui (\citeyear{CC07}), it is impractical to estimate $b$ by the method
of moments as it involves many terms and high-order moments. In
simulations, we used a robustified bootstrap estimate of $b$ suggested
by Chen and Cui (\citeyear{CC07}).

\section{Applications}

\subsection{Confidence regions for population mean}

A classical problem is the construction of confidence regions or
testing a hypothesis about a specific value of the population mean
based on a set of $n$ independent and identically distributed
observations. Particularly for scalar observations, the standard
approach is to use the Studentized sample mean,
\[
T_n(\theta) = \frac{\sqrt{n}( \bar x_n - \theta) }{s_n}
\]
for both purposes where $\bar x_n$ is the sample mean, and $s_n^2$ is
the sample variance. When the population distribution is normal,
$T_n(\theta)$ has a t-distribution with $n-1$ degrees of freedom. The
confidence interval or hypothesis test calibrated by the t-distribution
is found to be accurate even for nonnormal population distributions and
for moderate sample size $n$. For multivariate observations, the
$t$-statistic is replaced by Hotelling's $T^2$ defined as
\[
T_n^2 ( \theta) =n ( \bar{X}_n - \theta)^T S_n^{-1} ( \bar{X}_n -
\theta)
\]
with $\bar X_n$ the vector sample mean and $S_n$ the sample covariance
matrix. When the observations have a $p$-dimensional multivariate
normal distribution, $(n-p) T_n^2 ( \theta) /\{p(n-1)\}$ has an
$F$-distribution with $p$ and $n-p$ degrees of freedom. The
$F$-distribution often serves as a reference distribution for both
hypothesis tests and constructing confidence regions, whether or not
the normality assumption holds. Surprisingly, the normal-theory-based
confidence regions have reasonably accurate coverage probabilities even
when the sample sizes are small and the population distributions
deviate from the normal. Thus they serve as a good barometer to gauge
the performance of a new method.

%
%
\begin{table}
\caption{Bartlett correction factors
and their average estimates under linear regression models}\label{table4}
\begin{tabular*}{\tablewidth}{@{\extracolsep{\fill}}lcccccc@{}}
\hline
& \multicolumn{3}{c}{$\bolds{N(0,1)}$} & \multicolumn{3}{c@{}}{\textbf{Exp(1)}} \\
[-4pt]
& \multicolumn{3}{c}{\hrulefill} & \multicolumn{3}{c@{}}{\hrulefill} \\
$\bolds{n}$ & $\bolds{b}$ & $\bolds{b_n}$ & $\bolds{\tilde b_n}$ & $\bolds{b}$
& $\bolds{b_n}$ & $\bolds{\tilde b_n}$ \\
\hline
\phantom{0}30 & 3.55 & 2.39 & 3.56 & 7.98 & 2.61 & 5.39 \\
\phantom{0}50 & 3.53 & 2.74 & 3.61 & 7.92 & 3.35 & 6.16 \\
100 & 3.90 & 3.28 & 3.86 & 9.00 & 4.58 & 7.07 \\
\hline
\end{tabular*}
\end{table}

The EL and AEL counterparts
are obtained by letting $g(x; \theta) = x - \theta$.
For the sake of comparison, we use the same
simulation set-ups as in DiCiccio, Hall and Romano (\citeyear{DHR91}).
We investigate the coverage probabilities
of 90\%, 95\% and 99\% confidence intervals based on the
following methods:
\begin{enumerate}[(1)]
\item[(1)] Hotelling's $T^2$ (including the univariate case), T$^2$;
\item[(2)] The usual empirical likelihood, EL;
\item[(3)] Bartlett-corrected empirical likelihood with moment
estimate $b_n$, BEL;
\item[(4)] Adjusted empirical likelihood with moment estimate
$b_n$, AEL;
\item[(5)] Bartlett-corrected empirical likelihood with
$\tilde b_n$, BEL*;
\item[(6)] Adjusted empirical likelihood with $\tilde b_n$, AEL*;
\item[(7)] Bartlett-corrected empirical likelihood with
known $b$ value, BEL$_t$;
\item[(8)] Adjusted empirical likelihood with known $b$ value, AEL$_t$;
\item[(9)] Adjusted empirical likelihood with level of adjustment
$a_n = \frac{1}{2} \log n $, AEL$_0$.
\end{enumerate}

We generated 10,000 samples from four distributions: (a) the standard
normal; (b) an exponential distribution with mean 1; (c) a normal
mixture $0.2 N(5, 1) + 0.8 N(-1.25, 1)$; and (d) the $\chi_1^2$
distribution. The results are presented in Table \ref{table1} where
$0.2N_1 + 0.8N_2$ denotes the normal mixture distribution.

%
%
\begin{table}
\tabcolsep=0pt
\caption{Coverage probabilities for one-sample population mean}
\label{table1}
\begin{tabular*}{\tablewidth}{@{\extracolsep{\fill}}lccccccccccc@{}}
\hline
& $\bolds n$ & \textbf{Level} & $\bolds{T^2}$ & \textbf{EL} & \textbf{BEL}
& \textbf{AEL} & \textbf{BEL*} & \textbf{AEL*} & \textbf{BEL}$_{\bolds{t}}$
& \textbf{AEL}$_{\bolds{t}}$ & \textbf{AEL}$_{\bolds{0}}$ \\
\hline
$N(0,1)$ & 20 & 90 & 90.1 & 88.2 & 89.0 & 89.1 & 89.3 & \textbf{89.5} &
89.3 & 89.4 & 91.0 \\
& & 95 & 95.1 & 93.2 & 94.0 & 94.0 & 94.2 & \textbf{94.4} & 94.2 & 94.3 &
95.4 \\
& & 99 & 98.9 & 97.9 & 98.2 & 98.3 & 98.3 & \textbf{98.4} & 98.3 & 98.4 &
98.9 \\
& 30 & 90 & 90.2 & 89.0 & 89.7 & 89.8 & 90.0 & \textbf{90.0} & 89.9 & 89.9
& 91.1 \\
& & 95 & 95.5 & 94.3 & 94.9 & 94.9 & 95.0 & \textbf{95.0} & 95.0 & 95.0 &
95.8 \\
& & 99 & 99.1 & 98.7 & 98.8 & 98.8 & 98.8 & \textbf{98.8} & 98.9 & 98.9 &
99.1 \\
[4pt]
Exp(1) & 20 & 90 & 87.5 & 85.6 & 86.8 & 87.0 & 87.6 & \textbf{88.2} & 88.2
& 88.9 & 88.7 \\
& & 95 & 92.0 & 91.2 & 91.8 & 91.9 & 92.3 & \textbf{92.8} & 92.8 & 93.5 &
93.4 \\
& & 99 & 96.6 & 96.7 & 97.0 & 97.1 & 97.2 & \textbf{97.4} & 97.5 & 98.0 &
97.9 \\
& 30 & 90 & 87.6 & 86.7 & 87.7 & 87.8 & 88.2 & \textbf{88.5} & 88.6 & 88.9
& 89.0 \\
& & 95 & 92.8 & 92.3 & 92.9 & 93.0 & 93.3 & \textbf{93.6} & 93.7 & 93.9 &
94.0 \\
& & 99 & 97.1 & 97.6 & 97.9 & 97.9 & 98.0 & \textbf{98.0} & 98.2 & 98.3 &
98.4 \\
[4pt]
$0.2 N_1 + 0.8 N_2$
& 20 & 90 & 88.4 & 88.4 & 89.5 & 89.5 & 91.0 & \textbf{91.8} & 89.2 & 89.2
& 90.9 \\
& & 95 & 92.8 & 93.3 & 94.3 & 94.3 & 95.0 & \textbf{95.5} & 94.1 & 94.1 &
95.2 \\
& & 99 & 97.0 & 97.8 & 98.0 & 98.0 & 98.1 & \textbf{98.2} & 98.0 & 98.0 &
98.4 \\
& 30 & 90 & 88.7 & 89.1 & 89.9 & 89.9 & 90.3 & \textbf{90.4} & 89.7 & 89.8
& 91.2 \\
& & 95 & 93.7 & 94.4 & 94.9 & 94.9 & 95.3 & \textbf{95.5} & 94.7 & 94.7 &
95.6 \\
& & 99 & 97.8 & 98.8 & 99.1 & 99.1 & 99.2 & \textbf{99.3} & 99.0 & 99.0 &
99.3 \\
[4pt]
$\chi^2_1$ & 20 & 90 & 84.8 & 83.7 & 85.0 & 85.2 & 86.4 & \textbf{87.3} &
87.2 & 89.2 & 86.7 \\
& & 95 & 89.2 & 89.3 & 90.4 & 90.5 & 91.3 & \textbf{92.0} & 92.2 & 93.8 &
91.7 \\
& & 99 & 94.4 & 95.4 & 96.0 & 96.0 & 96.4 & \textbf{96.8} & 96.9 & 98.5 &
96.9 \\
& 30 & 90 & 85.9 & 85.4 & 86.5 & 86.7 & 87.7 & \textbf{88.2} & 88.2 & 88.9
& 87.8 \\
& & 95 & 90.2 & 91.1 & 91.9 & 91.9 & 92.4 & \textbf{92.7} & 93.0 & 93.6 &
92.8 \\
& & 99 & 95.2 & 96.5 & 96.8 & 96.8 & 97.0 & \textbf{97.2} & 97.3 & 97.7 &
97.3 \\
\hline
\end{tabular*}
\end{table}

Under the normal model, $T^2$ is optimal, yet we find
that the AEL* is as good within simulation error.
The accuracy of the AEL* is consistently better than
that of the BEL and BEL*.
This is particularly true when the population distribution
is exponential or chi-square. Under the mixture model,
the AEL* has a slightly higher than nominal coverage
probability.
Finally, we remark that under the chi-square distribution,
all the methods still have room for improvement when
$n=20$. Our simulation results on EL and BEL
are comparable
to those reported in the literature.

In the multivariate case, we conducted simulation experiments
for $p=2$ and $p=3$. We used the following strategy
to generate correlated trivariate observations.
We first generated a random observation $D$ from
the uniform distribution on the interval [1, 2].
Given $D$, we generated $X_1, X_2$ and $X_3$
from the distributions specified as follows:

\begin{enumerate}[(a)]
\item[(a)] $X_1 \sim N(0,D^2)$,
$X_2 \sim \operatorname{Gamma}(D^{-1},1)$, $X_3 \sim\chi_{D}^2$;
\item[(b)] $X_1 \sim \operatorname{Gamma}(D,1)$,
$X_2 \sim \operatorname{Gamma}(D^{-1},1)$,
$X_3 \sim \operatorname{Gamma}(4-D,1)$;
\item[(c)] $X_1 \sim0.2N(5,D^2)+0.8N(-1.25, D^{-2})$,
$X_2 \sim0.2N(5,D^{-2})+0.8N(-1.25,\break D^2)$,
$X_3 \sim N(0, D^2)$;
\item[(d)] $X_1 \sim \operatorname{Poisson}(D)$,
$X_2 \sim \operatorname{Poisson}(D^{-1})$,
$X_3 \sim \operatorname{Poisson}(4-D)$.
\end{enumerate}

When $p=2$, we used $X_1$ and $X_2$ in our simulation and generated
10,000 data sets with sample sizes $n=20$ and 30. When $p=3$, we also
generated 10,000 data sets but increased sample sizes to $n=30$ and
$50$ to accommodate the higher dimension. Table \ref{multi1} presents
the simulation results.

%
%
\begin{table}
\tabcolsep=0pt
\caption{Coverage probabilities for one-sample
multivariate ($q=2, 3$) population mean}\label{multi1}
\begin{tabular*}{\tablewidth}{@{\extracolsep{\fill}}lcccccccccccc@{}}
\hline
&& $\bolds n$ & \textbf{Level} & $\bolds{T^2}$ & \textbf{EL} & \textbf{BEL}
& \textbf{AEL} & \textbf{BEL*} & \textbf{AEL*} & \textbf{BEL}$_{\bolds{t}}$
& \textbf{AEL}$_{\bolds{t}}$ & \textbf{AEL}$_{\bolds{0}}$ \\
\hline
$q=2$& (a) & 20 & 90 & 86.0 & 81.7 & 83.8 & 84.3 & 84.8 & \textbf{86.2} &
85.4 & 87.0 & 86.6 \\
&& & 95 & 91.3 & 87.7 & 89.3 & 89.8 & 90.1 & \textbf{91.6} & 90.6 &
92.2 &
91.8 \\
&& & 99 & 96.5 & 94.5 & 95.3 & 95.9 & 95.9 & \textbf{96.7} & 96.1 &
97.9 &
97.4 \\
&& 30 & 90 & 87.2 & 85.1 & 86.5 & 86.7 & 87.0 & \textbf{87.8} & 87.4 &
88.0 & 88.2 \\
&& & 95 & 92.2 & 90.8 & 91.7 & 92.0 & 92.2 & \textbf{92.8} & 92.4 &
93.0 &
93.2 \\
&& & 99 & 97.0 & 96.5 & 97.0 & 97.2 & 97.2 & \textbf{97.5} & 97.4 &
97.6 &
97.8 \\
[4pt]
&(b) & 20 & 90 & 84.5 & 80.8 & 82.7 & 83.4 & 84.2 & \textbf{86.2} &
84.9 &
87.2 & 85.6 \\
&& & 95 & 89.5 & 86.8 & 88.4 & 89.1 & 89.6 & \textbf{91.1} & 90.0 &
92.5 &
91.1 \\
&& & 99 & 95.4 & 93.6 & 94.5 & 95.0 & 94.9 & \textbf{96.2} & 95.3 &
98.5 &
96.7 \\
&& 30 & 90 & 85.9 & 84.5 & 86.0 & 86.3 & 86.8 & \textbf{87.6} & 87.1 &
88.0 & 87.6 \\
&& & 95 & 90.7 & 90.4 & 91.6 & 91.8 & 92.2 & \textbf{92.7} & 92.5 &
93.1 &
92.9 \\
&& & 99 & 96.1 & 96.3 & 96.8 & 97.0 & 97.1 & \textbf{97.4} & 97.3 &
97.8 &
97.6 \\
[4pt]
&(c) & 20 & 90 & 85.7 & 84.6 & 86.2 & 86.4 & 87.7 & \textbf{89.4} &
86.2 &
86.5 & 88.8 \\
&& & 95 & 90.6 & 89.9 & 91.1 & 91.5 & 92.3 & \textbf{93.7} & 91.2 &
91.4 &
93.2 \\
&& & 99 & 95.8 & 95.2 & 95.7 & 96.0 & 96.0 & \textbf{97.2} & 95.7 &
96.0 &
97.1 \\
&& 30 & 90 & 87.9 & 87.4 & 88.9 & 89.0 & 89.6 & \textbf{90.0} & 88.9 &
89.0 & 90.7 \\
&& & 95 & 92.9 & 93.2 & 94.2 & 94.3 & 94.7 & \textbf{95.1} & 94.1 &
94.2 &
95.4 \\
&& & 99 & 97.2 & 98.0 & 98.3 & 98.4 & 98.5 & \textbf{98.8} & 98.3 &
98.4 &
98.7 \\
[4pt]
&(d) & 20 & 90 & 88.5 & 84.2 & 85.9 & 86.2 & 86.6 & \textbf{87.4} &
86.8 &
87.6 & 89.0\\
&& & 95 & 93.3 & 90.2 & 91.3 & 91.6 & 91.8 & \textbf{92.5} & 91.8 &
92.4 &
93.7 \\
&& & 99 & 97.6 & 95.8 & 96.2 & 96.5 & 96.4 & \textbf{97.0} & 96.5 &
97.0 &
98.1 \\
&& 30 & 90 & 88.4 & 86.4 & 87.6 & 87.7 & 87.9 & \textbf{88.2} & 88.0 &
88.3 & 89.8 \\
&& & 95 & 93.6 & 92.3 & 93.0 & 93.1 & 93.3 & \textbf{93.5} & 93.3 &
93.5 &
94.2 \\
&& & 99 & 98.0 & 97.2 & 97.6 & 97.7 & 97.8 & \textbf{97.9} & 97.8 &
97.9 &
98.4 \\
\hline
\end{tabular*}
\end{table}

\setcounter{table}{4}
\begin{table}
\tabcolsep=0pt
\caption{(Continued)}
\begin{tabular*}{\tablewidth}{@{\extracolsep{\fill}}lcccccccccccc@{}}
\hline
&& $\bolds n$ & \textbf{Level} & $\bolds{T^2}$ & \textbf{EL} & \textbf{BEL}
& \textbf{AEL} & \textbf{BEL*} & \textbf{AEL*} & \textbf{BEL}$_{\bolds{t}}$
& \textbf{AEL}$_{\bolds{t}}$ & \textbf{AEL}$_{\bolds{0}}$ \\
\hline
$q=3$&(a) & 30 & 90 & 85.2 & 81.5 & 83.8 & 84.6 & 84.9 & \textbf{86.5} &
85.3 & 87.2 & 86.0 \\
&& & 95 & 90.6 & 88.1 & 89.8 & 90.7 & 90.6 & \textbf{91.9} & 91.0 &
92.5 &
91.6 \\
&& & 99 & 96.2 & 94.8 & 95.6 & 96.4 & 96.1 & \textbf{97.1} & 96.2 &
97.9 &
97.0 \\
&& 50 & 90 & 85.8 & 84.4 & 85.8 & 86.2 & 86.5 & \textbf{86.9} & 86.5 &
87.0 & 86.9 \\
&& & 95 & 91.2 & 90.7 & 91.9 & 92.2 & 92.2 & \textbf{92.7} & 92.4 &
92.8 &
92.7 \\
&& & 99 & 96.6 & 96.6 & 97.2 & 97.5 & 97.5 & \textbf{97.7} & 97.5 &
97.8 &
97.7 \\
[4pt]
&(b) & 30 & 90 & 85.3 & 81.4 & 83.6 & 84.4 & 84.8 & \textbf{86.1} &
85.2 &
86.7 & 86.0 \\
&& & 95 & 90.8 & 87.8 & 89.7 & 90.3 & 90.4 & \textbf{91.7} & 90.8 &
92.3 &
91.6 \\
&& & 99 & 96.4 & 95.1 & 95.9 & 96.5 & 96.3 & \textbf{97.1} & 96.4 &
97.6 &
97.2 \\
&& 50 & 90 & 86.7 & 85.7 & 87.1 & 87.4 & 87.6 & \textbf{88.0} & 87.7 &
88.1 & 88.1 \\
&& & 95 & 92.0 & 91.1 & 92.2 & 92.5 & 92.6 & \textbf{92.8} & 92.8 &
93.1 &
93.1 \\
&& & 99 & 97.1 & 97.5 & 97.8 & 97.9 & 97.9 & \textbf{98.0} & 98.0 &
98.1 &
98.2 \\
[4pt]
&(c) & 30 & 90 & 88.0 & 84.7 & 86.7 & 87.0 & 87.2 & \textbf{88.0} &
86.9 &
87.4 & 88.8 \\
&& & 95 & 93.0 & 90.5 & 91.9 & 92.3 & 92.4 & \textbf{93.1} & 92.1 &
92.5 &
93.7 \\
&& & 99 & 97.6 & 96.5 & 97.0 & 97.3 & 97.2 & \textbf{98.0} & 97.1 &
97.3 &
98.1 \\
&& 50 & 90 & 88.7 & 87.4 & 88.7 & 88.8 & 89.0 & \textbf{89.1} & 88.8 &
88.9 & 90.0 \\
&& & 95 & 93.5 & 93.2 & 94.1 & 94.2 & 94.3 & \textbf{94.4} & 94.2 &
94.2 &
94.9 \\
&& & 99 & 98.2 & 98.3 & 98.6 & 98.6 & 98.6 & \textbf{98.7} & 98.6 &
98.6 &
98.9 \\
[4pt]
&(d) & 30 & 90 & 88.4 & 84.2 & 86.1 & 86.6 & 86.7 & \textbf{87.3} &
86.7 &
87.3 & 88.4 \\
&& & 95 & 93.7 & 90.5 & 91.9 & 92.3 & 92.3 & \textbf{93.0} & 92.4 &
93.0 &
93.7 \\
&& & 99 & 98.1 & 96.4 & 97.2 & 97.4 & 97.3 & \textbf{97.7} & 97.3 &
97.7 &
98.3 \\
&& 50 & 90 & 88.7 & 86.8 & 88.2 & 88.3 & 88.4 & \textbf{88.5} & 88.4 &
88.5 & 89.4 \\
&& & 95 & 94.0 & 92.9 & 93.7 & 93.8 & 93.8 & \textbf{93.9} & 93.8 &
93.9 &
94.4 \\
&& & 99 & 98.4 & 97.9 & 98.2 & 98.3 & 98.3 & \textbf{98.3} & 98.3 &
98.3 &
98.6 \\
\hline
\end{tabular*}
\end{table}

We observe that the AEL* outperforms all other methods, often
substantially. Under the bivariate mixture model (c) at nominal level
95\% and sample size $n=20$, the AEL* has 93.7\% coverage probability
compared to 91.1\% for the BEL and 92.3\% for the BEL*. This is
significant because the AEL*, the BEL and the BEL* are known to be
precise up to the same order $n^{-2}$. The difference in performances
presumably comes from higher orders.

We remark here that the above discussion has not taken AEL$_t$ and
AEL$_0$ into account. The AEL$_t$ is only for theoretical interest and
its performance indicates how far AEL can be improved by choosing a
better estimator of $b$. The AEL$_0$ is the AEL with a conventional
level of adjustment suggested in Chen, Variyath and Abraham
(\citeyear{CVA08}). It has comparable performance to AEL*. Due to a
lack of theoretical justification, the observed good performance is
hard to generalize. We will continue to keep an eye on its performance.

\subsection{Linear regression}

The empirical likelihood method can also be used to construct
confidence regions for the regression coefficient $\beta$ in the
following linear regression model:
%
%
\begin{equation}
\label{reg}
y = \mathbf{x}^T \beta+\varepsilon,
\end{equation}
where $\beta$ is a $p$-dimensional parameter, $\mathbf{x}$ a
$p$-dimensional fixed design point and $y$ the scalar response. Chen
(\citeyear{C93}) showed that the empirical likelihood confidence
regions for $\beta$ are also Bartlett correctable. In comparison, by
letting $g(y, \mathbf{x}; \beta) = \mathbf{x}( y - \mathbf{x}^T \beta)$
the proposed AEL method (AEL*) directly applies.

%
%
\begin{table}
\tabcolsep=0pt
\caption{Coverage probabilities for the regression coefficient $\beta$}
\label{table3}
\begin{tabular*}{\tablewidth}{@{\extracolsep{\fill}}lcccccccccd{3.1}c@{}}
\hline
& $\bolds n$ & \textbf{Level} & $\bolds F$\textbf{-test} & \textbf{EL}
& \textbf{BEL} & \textbf{AEL} & \textbf{BEL*} & \textbf{AEL*}
& \textbf{BEL}$_{\bolds{t}}$ & \multicolumn{1}{c}{\textbf{AEL}$_{\bolds{t}}$}
& \textbf{AEL}$_{\mathbf{0}}$ \\
\hline
$N(0,1)$
& \phantom{0}30 & 90 & 90.0 & 84.0 & 85.7 & 86.1 & 86.6 & \textbf{87.7} & 86.6 & 87.5
& 87.4\\
& & 95 & 94.9 & 90.1 & 91.5 & 92.0 & 92.2 & \textbf{93.3} & 92.3 & 93.2 &
93.0\\
& & 99 & 99.3 & 96.6 & 97.3 & 97.5 & 97.4 & \textbf{98.2} & 97.5 & 98.2 &
98.1\\
& \phantom{0}50 & 90 & 89.7 & 86.9 & 88.4 & 88.5 & 88.7 & \textbf{88.9} & 88.7 & 89.0
& 89.2\\
& & 95 & 95.0 & 92.7 & 93.6 & 93.7 & 93.8 & \textbf{94.0} & 93.8 & 94.0 &
94.2\\
& & 99 & 99.0 & 97.7 & 98.1 & 98.2 & 98.2 & \textbf{98.2} & 98.2 & 98.3 &
98.4\\
&100 & 90 & 89.6 & 88.3 & 89.1 & 89.1 & 89.2 & \textbf{89.2} & 89.2 & 89.3
& 89.4\\
& & 95 & 94.8 & 93.8 & 94.3 & 94.4 & 94.4 & \textbf{94.5} & 94.4 & 94.5 &
94.5\\
& & 99 & 99.0 & 98.5 & 98.6 & 98.7 & 98.7 & \textbf{98.7} & 98.7 & 98.7 &
98.8\\
[4pt]
Exp(1) & \phantom{0}30 & 90 & 87.9 & 79.6 & 81.9 & 82.4 & 83.6 & \textbf{86.1} & 85.7
& 92.6 & 83.5\\
& & 95 & 92.8 & 86.4 & 88.2 & 88.8 & 89.4 & \textbf{91.6} & 91.0 & 98.5 &
89.6\\
& & 99 & 97.7 & 93.7 & 94.7 & 95.2 & 95.3 & \textbf{97.0} & 96.3 &
100.0 &
95.9\\
& \phantom{0}50 & 90 & 88.7 & 83.7 & 85.4 & 85.6 & 86.4 & \textbf{87.5} & 87.4 & 89.0
& 86.0 \\
& & 95 & 93.8 & 90.0 & 91.3 & 91.5 & 92.1 & \textbf{92.8} & 92.9 & 94.2 &
91.8 \\
& & 99 & 98.3 & 96.3 & 96.9 & 97.1 & 97.3 & \textbf{97.8} & 97.7 & 98.9 &
97.2 \\
&100 & 90 & 88.9 & 86.2 & 87.3 & 87.3 & 87.8 & \textbf{88.1} & 88.4 & 88.8
& 87.4\\
& & 95 & 94.2 & 92.2 & 93.0 & 93.0 & 93.3 & \textbf{93.6} & 93.8 & 94.2 &
93.1\\
& & 99 & 98.5 & 97.8 & 98.1 & 98.1 & 98.2 & \textbf{98.3} & 98.3 & 98.5 &
98.1\\
\hline
\end{tabular*}
\end{table}

In this simulation study, we examined the performance of the AEL*
method based on model (\ref{reg}) with $p=2$, the true parameter value
$\beta_0=(1,1)^T$, and the errors $\varepsilon_i$ were generated from
either a normal distribution or from a centralized exponential
distribution as specified in Table \ref{table3}. The design matrix of
$\mathbf{x}$ of size $n \times2$ was taken from the first $n$ rows in
Table 1 of Chen (\citeyear{C93}). The simulation results also are given
in Table \ref{table3}. The improvement of the AEL* over the EL, BEL or
BEL* is universal and substantial, particularly under the nonnormal
models when the sample sizes are small.

\subsection{An example where $q > p$}

In this subsection, we examine the AEL through an asset-pricing model
investigated by Hall and Horowitz (\citeyear{HH96}) and also by Imbens,
Spady and Johnson (\citeyear{ISJ98}) expanded with $q$ $(q \geq2)$
moment restrictions by Schennach (\citeyear{S07}). The parameter of
interest is defined through the following estimating equations:
%
%
\begin{equation}
\label{gee} E g(X; \theta) \equiv E \pmatrix{
r(X,\theta) \cr
X_2 r(X,\theta) \cr
( X_3 -1 ) r(X,\theta) \cr
\vdots\cr
( X_q -1 )r(X,\theta)}
=0,
\end{equation}
where $ r(X,\theta) = \exp\{ -0.72 - \theta( X_1 + X_2 ) + 3 X_2 \} -
1, $ $X=(X_1,X_2,\ldots, X_q)$ and $\theta$ is a scalar parameter.
Components of $X$ are mutually independent and $X_1, X_2
\stackrel{\mathrm{i.i.d.}}{\sim} N(0,0.16)$, $X_3, \ldots, X_q
\stackrel{\mathrm{i.i.d.}}{\sim} \chi_1^2$. We generated data from the
models with $\theta_0 = 3$, $q=2$ and $q=3$, respectively.

Although Theorem \ref{thm2} is applicable, precisely estimating $b$ is
not easy due to its complex expression. Instead, Chen and Cui
(\citeyear{CC07}) proposed a bootstrap estimate. We adopted\vspace*{2pt} their
strategy with a robust modification. Let $\Delta_m$ be the sample
median of $\Delta_n^*(\hat\theta; 0)$ based on $B=300$ bootstrap
samples. We estimate $b$ by
\[
\hat b = n ( \Delta_m / 0.4549 - 1),
\]
where $0.4549$ is the median of the $\chi_1^2$ distribution. We
generated samples of sizes $n= 100 $ and $200$. The average bootstrap
estimates of $\hat b$ are $31$ and $58$ for $q = 2$ and $q=3$ over
$1000$ repetitions. We call them off-line estimates of $b$ and carried
out the corresponding simulations side-by-side with the bootstrap
estimator $\hat b$ for each sample generated.

In Table \ref{table5}, we report the coverage probabilities of the
nominal 90\%, 95\% and 99\% confidence intervals of the empirical
likelihood (EL), the Bartlett corrected empirical likelihood (BEL), the
adjusted empirical likelihood [AEL(5)] and the adjusted empirical
likelihood with conventional $a_n = \log(n)/2$ (AEL$_0$). Due to the
exponential nature of $g$ in $\theta$ in this example, the sample mean
$\bar g$ is unstable. For robustness, we computed $g_{n+1}$ with the
trimmed mean by removing five largest $\|g_i\|$ values.

%
%
\begin{table}
\caption{Simulation results under the expanded asset-pricing model}
\label{table5}
\begin{tabular*}{\tablewidth}{@{\extracolsep{\fill}}lcccccccc@{}}
\hline
 & \textbf{Level} & \textbf{EL} & \textbf{BEL} &
\textbf{AEL(5)} & & \textbf{BEL} & \textbf{AEL(5)} & \textbf{AEL}$_{\mathbf{0}}$
\\
\hline
& && \multicolumn{2}{c}{\textbf{Bootstrapped} $\bolds b$}& & \multicolumn{2}{c}{\textbf{Off-line} $\bolds{b=31}$} \\[-4pt]
& && \multicolumn{2}{c}{\hrulefill}& & \multicolumn{2}{c}{\hrulefill} \\
$q=2$\\
\quad$n=100$ & 90 & 82.6 & 86.5 & 85.3 & & 87.4 & 89.8 & 82.7 \\
& 95 & 88.4 & 91.2 & 92.6 & & 92.8 & 95.4 & 88.8 \\
& 99 & 95.8 & 96.7 & 97.3 & & 97.3 & 99.5 & 95.9 \\
[4pt]
\quad$n=200$ & 90 & 83.9 & 86.6 & 85.1 & & 87.8 & 87.2 & 84.3 \\
& 95 & 91.2 & 92.6 & 91.9 & & 93.1 & 93.3 & 91.4 \\
& 99 & 96.9 & 97.4 & 97.6 & & 97.8 & 98.2 & 96.9 \\
[4pt]
& && \multicolumn{2}{c}{\textbf{Bootstrapped} $\bolds b$}& & \multicolumn{2}{c}{\textbf{Off-line} $\bolds{b=58}$} \\[-4pt]
& && \multicolumn{2}{c}{\hrulefill}& & \multicolumn{2}{c}{\hrulefill} \\
$q=3$\\
\quad$n=100$ & 90 & 78.4 & 84.9 & 84.1 & & 87.4 & 90.5 & 79.8 \\
& 95 & 85.7 & 90.8 & 90.4 & & 93.1 & 96.7 & 86.1 \\
& 99 & 94.0 & 96.1 & 97.9 & & 97.7 & 99.8 & 94.0 \\
[4pt]
\quad$n=200$ & 90 & 82.5 & 86.9 & 86.5 & & 87.4 & 89.8 & 82.5 \\
& 95 & 89.7 & 92.7 & 92.9 & & 93.3 & 95.3 & 89.8 \\
& 99 & 96.1 & 97.2 & 98.5 & & 97.6 & 99.2 & 96.1 \\
\hline
\end{tabular*}
\end{table}

In terms of the precision of the coverage probabilities, the AEL is
better than the BEL which is better than the EL and the AEL$_0$, and
the latter two have similar performances. Even after the
robustification, the bootstrap estimation of $b$ ranges from $-27$ to
$376$ when $n=100$. This observation indicates that neither the BEL nor
the AEL is ready to be applied to models similar to the one in this
example. The simulation results have instead shown the potential of the
AEL approach. We hope to further investigate this problem in the
future.

\begin{appendix}\label{app}
\section*{Appendix}
\begin{pf*}{Proof of Theorem \protect\ref{thm1}}
We now present the proof for the general
case where $g(x; \theta)$ is vector valued.

In addition to the notation introduced earlier, we further define
\[
A^{rs\cdots t}
=
\frac{1}{n}\sum_i^n Y^{r}Y^{s}\cdots Y^{t} - \alpha^{rs\cdots t},
\]
where $ \alpha^{rs\cdots t} $ is defined in (\ref{alpha}). Without loss
of generality, we assume that $\alpha^{rs} = I(r = s)$ at $\theta=
\theta_0$. By DiCiccio, Hall and Romano (\citeyear{DHR91}), the
solution to (\ref{eqn4}), before any adjustment, can be expanded as
\[
\lambda
=
\lambda_1 + \lambda_2 + \lambda_3 + O_p(n^{-2})
\]
with
\[
\lambda_1^r = A^r,\qquad \lambda_2^r = -A^{rs} A^s + \alpha^{rst} A^s A^t
\]
and
\begin{eqnarray*}
\lambda_3^r &=& A^{rs} A^{tu} A^u + A^{rst} A^s A^t
+ 2 \alpha^{rst} \alpha^{tuv} A^s A^u A^v \\
&&{} - 3 \alpha^{rst} A^{tu} A^s A^u - \alpha^{rstu} A^s A^t A^u.
\end{eqnarray*}
Here we have used the summation convention according to which, if an
index occurs more than once in an expression, summation over the index
is understood. Substituting these expansions into the expression for
$R_n(\theta_0)$, we get
%
%
\begin{equation}
R_n(\theta_0) =n \{ R_1 + R_2 + R_3\}^T \{ R_1 + R_2 + R_3\} +
O_p(n^{-3/2})
\end{equation}
with
%
%
\begin{eqnarray}
\label{R1}
R_1^r &=& A^r,\qquad
R_2^r = \tfrac{1}{3} \alpha^{rst} A^s A^t - \tfrac{1}{2} A^{rs}A^s,
\nonumber\\
R_3^r &=& \tfrac{3}{8} A^{rs} A^{st} A^t - \tfrac{5}{12} \alpha^{rst}
A^{tu} A^s A^u - \tfrac{5}{12} \alpha^{stu} A^{rs} A^t A^u
\\
&&{} + \tfrac{4}{9} \alpha^{rst} \alpha^{tuv} A^s A^u A^v +
\tfrac{1}{3} A^{rst} A^s A^t - \tfrac{1}{4} \alpha^{rstu} A^s A^t
A^u.\nonumber
\end{eqnarray}

Recall the usual Lagrange multiplier $\lambda$ solves
$f(\lambda) = 0$ where
\[
f(\zeta) = n^{-1}\sum_{i=1}^n \frac{g(x_i; \theta)}{1+ \zeta^Tg(x_i;
\theta)}.
\]
Now we work on the Lagrange multiplier after an adjustment at level
$a_n = a + O_p(n^{-1/2})$. Since $\lambda_a = O_p(n^{-1/2})$, it must
solve
\[
f(\lambda_a) - \frac{a}{n} \bar g = O_p(n^{-2}).
\]
A Taylor expansion of $f(\lambda_a)$ gives
\[
f(\lambda_a) = f(\lambda) + \frac{\partial f(\lambda)}{\partial\lambda}
(\lambda _a - \lambda) + O\bigl((\lambda_a - \lambda)^2\bigr).
\]
Since $f(\lambda)=0$, it simplifies to
\[
\lambda_a - \lambda = \frac{a}{n} \biggl( \frac{\partial
f(\lambda)}{\partial\lambda} \biggr)^{-1} \bar g + O_p(n^{-2}).
\]
Note that
\[
\frac{\partial f(\lambda)}{\partial\lambda} = - E \{ g(X;
\theta_0)g^T(X; \theta_0) \} + O_p(n^{-1/2})
\]
and by assumption $E \{ g(X; \theta_0)g^T(X; \theta_0) \} =I $; thus we
arrive at
\[
\lambda_a = \lambda- {n^{-1} a} \bar g + O_p(n^{-2}) = ( 1 - {n^{-1} a}
) \lambda+ O_p(n^{-2}).
\]
That is, the two Lagrange multipliers are nearly equal.

Next, we quantify the effect of slightly different Lagrange multipliers
on the expansion of $R_n(\theta_0; a_n)$. We have
\begin{eqnarray*}
R_n(\theta_0; a_n)
&=&
2 \sum_{i=1}^n \log\{ 1 +(1 - n^{-1} a) \lambda^T g_i\} \\
&&{} + 2 \log\{ 1 - (1 - n^{-1} a) a \lambda^T \bar g\} + O(n^{-3/2}).
\end{eqnarray*}
Note that
\[
\log\{ 1 - (1 - n^{-1} a) a \lambda^T \bar g\} = - a \lambda^T \bar g +
O_p(n^{-2})
\]
and, surprisingly,
\[
2 \sum_{i=1}^n \log\{ 1 +(1 - n^{-1} a) \lambda^T g_i \} =
R_n(\theta_0) + O_p(n^{-3/2}).
\]
Therefore, we must have
%
%
\begin{equation}
R(\theta_0; a_n) = R_n(\theta_0) - 2 a R_1^T R_1 + O_p(n^{-3/2})
\end{equation}
where $R_1$ is defined in (\ref{R1}), and, consequently,
\[
R_n(\theta_0; a_n) = n \{R_1+R_2 + R_{3a}\}^T \{R_1+R_2 + R_{3a}\} +
O_p(n^{-3/2})
\]
with
%
%
\begin{equation}
\label{Ra}
R_{3a} = R_3 - n^{-1} a R_1.
\end{equation}

Denote
\begin{eqnarray*}
Q_n &=&
\sqrt{n}(R_1+R_2 + R_{3a}), \\
U_n &=& (A^1, \ldots, A^q, A^{11}, A^{12}, \ldots, A^{qq}, A^{111},
A^{112}, \ldots, A^{qqq} )^T
\end{eqnarray*}
such that the super-indices in $A^{rst}$ satisfy $1 \leq r \leq s \leq
t \leq q$. Hence, $U_n$ has $q(q+1)(q+2)/6$ components, and each
component is a centralized sample mean. Furthermore, $Q_n$ is a smooth
vector-valued function of $U_n$. According to Bhattacharya and Ghosh
(\citeyear{BG78}), the Edgeworth expansion of a smooth function of the
sample mean (vector valued) is given by its formal Edgeworth expansion
based on its cumulants. Depending on the required order of the
expansion, the appropriate lower-order cumulants must exist.

In this theorem, we look for an expansion of the density function of
$Q_n$ up to order $o(n^{-2})$. This expansion is determined by the
first six cumulants of $U_n$ and the derivative of $Q_n$ with respect
to $U_n$. Note that we assumed that the 18th moment of $g(x; \theta)$
exists and the highest order in $U_n$ is three, hence all cumulants of
$U_n$ up to order $6$ exist. The cumulants of $Q_n$ can then be
obtained through those of $U_n$.

Let $\kappa_{r, s, \ldots, t}(Q_n)$ denote the joint cumulant of the
$r$th, $s$th$, \ldots, t$th components of $Q_n$. After some lengthy but
routine algebraic work, we get
\begin{eqnarray*}
\kappa_{r}(Q_n)
&=& - n^{-1/2} \mu^r + n^{-3/2} c_1^r + o(n^{-2}), \\
\kappa_{r, s}(Q_n)
&=& I(r = s) + n^{-1} \gamma^{rs} + n^{-2} c_2^{rs}+o(n^{-2}), \\
\kappa_{r, s, t}(Q_n)
&=& n^{-3/2} c_3^{rst} + o(n^{-2}), \\
\kappa_{r, s, t, u}(Q_n)
&=& n^{-2} c_4^{rstu} + o(n^{-2}),
\end{eqnarray*}
where
\begin{eqnarray*}
\mu^r &=&
\tfrac{1}{6} \alpha^{rss}, \\
\gamma^{rs} &=& \tfrac{1}{2} \alpha^{rstt} -
\tfrac{1}{3}\alpha^{rtu}\alpha^{stu} - \tfrac{1}{36}
\alpha^{rst}\alpha^{tuu} - 2 a I(r=s)
\end{eqnarray*}
and $c_1^r$, $c_2^{rs}$, $c_3^{rst}$, $c_4^{rstu}$ are some nonrandom
constants. Cumulants of orders five and six are $o(n^{-2})$.

Let $f_{Q_n}(\mathbf{z})$ and $\phi(\mathbf{z})$ be the density
functions of $Q_n$ and the $q$-variate standard normal distribution.
The key consequence of the above computation is the resultant formal
Edgeworth expansion,
\[
f_{Q_n}(\mathbf{z}) = \Biggl\{1 + \sum_{i=1}^4 n^{-i/2} \pi_i(\mathbf{z})+
o(n^{-2}) \Biggr\} \phi ({\mathbf{z}})
\]
with
\begin{eqnarray*}
\pi_1(\mathbf{z})
&=& \mu^r \mathbf{z}^r, \\
\pi_2(\mathbf{z})
&=& \tfrac{1}{2}(\gamma^{rs} + \mu^r\mu^s)\{\mathbf{z}^r\mathbf
{z}^s - I(r=s)\}
\end{eqnarray*}
and for some polynomials $\pi_3(\mathbf{z})$ and $\pi_4(\mathbf{z})$
which are of order no more than four, the former is odd and the latter
is even. Their specific forms are not needed further and so are
omitted.

The above expansion implies that
\[
\mbox{\textsc{pr}}\{ Q_n^T Q_n \leq x \} = \int_{\mathbf{z}^T
\mathbf{z}< x} \Biggl\{ 1 + \sum_{i=1}^4 n^{-i/2} \pi_i(\mathbf{z}) \Biggr\}
\phi(\mathbf{z}) \, d\mathbf{z}+ o(n^{-2}).
\]
Because $\pi_1(\mathbf{z})$ and $\pi_3(\mathbf{z})$ are odd functions,
their integrations over the symmetric region are zero. For the same
reason, the integrations of the $\mathbf{z}^r \mathbf{z}^s$ terms in
$\pi_2(\mathbf{z})$ when $r \neq s$ over a symmetric region are also
zero. We further note that the expression of $\gamma^{rs}$ involves
$a$, and it is simple to get
\[
\int_{\mathbf{z}^T \mathbf{z}< x} \pi_2(\mathbf{z}) \phi({\mathbf {z}})
\,d\mathbf{z} = \frac{1}{2} ( b - 2 a ) \int_{\mathbf{z}^T \mathbf{z}< x}
(\mathbf{z}^T \mathbf{z}- q) \phi(\mathbf{z}) \,d\mathbf{z},
\]
where
\[
b = \frac{1}{q} \biggl( \frac{1}{2}\alpha^{rrss} - \frac{1}{3}\alpha
^{rst}\alpha^{rst}\biggr).
\]
This $b$ is the Bartlett correction factor given in DiCiccio, Hall and
Romano (\citeyear{DHR91}). Its expression is simpler than the earlier
one because we assumed $\alpha^{rs} = I(r = s)$. Hence, when $ a =
b/2$, we have
\[
\mbox{\textsc{pr}}\{ Q_n^T Q_n \leq x \} = \int_{\mathbf{z}^T
\mathbf{z}< x} \phi(\mathbf{z}) \, d\mathbf{z}+ O(n^{-2}) =
\mbox{\textsc{pr}}(\chi_q^2 \leq x) + O(n^{-2}).
\]
This completes the proof.
\end{pf*}

The conclusion for $R_n(\theta_0; a_{1n}, a_{2n})$ is obtained
similarly.
\begin{pf*}{Proof of Theorem \protect\ref{thm2}}
Expanding $\Delta_n(\theta_0; a_n)$ and then computing its cumulants
are by far the most demanding parts of the proof of Theorem \ref{thm2}.
The tasks are formidable. Fortunately, we find a short-cut by relating
$\Delta_n(\theta_0; a_n)$ to $\Delta_n(\theta_0; 0)$. By Chen and Cui
(\citeyear{CC07}),
\begin{eqnarray*}
\Delta_n(\theta_0; 0) & = &
R_n(\theta_0; 0) - \inf_\theta R_n(\theta; 0)\\
&= & n \{R_1+R_2 + R_3\}^T \{R_1+R_2 + R_3\} + O_p(n^{-3/2})
\end{eqnarray*}
for some $R_1, R_2 $ and $R_3 $; some of which are different from those
in DiCiccio, Hall and Romano (\citeyear{DHR91}). They have the same
fundamental properties that enable the Bartlett correction. In
addition, $R_1$ equals the first $p$ components of $n^{-1} \sum_{i=1}^n
g(X_i; \theta_0)$ after $g$ is standardized in some way.

With some relatively routine algebra, we find
\[
R_n(\theta_0; a_n) = R_n(\theta_0; 0) - 2a \sum_{r=1}^{q} \Biggl\{ n^{-1}
\sum_{i=1}^n g^r(X_i; \theta_0)\Biggr\}^2 + O_p(n^{-3/2})
\]
and
\[
\inf_\theta R_n(\theta; a_n) = \inf_\theta R_n(\theta; 0) - 2a \sum_{
r=p+1}^{ q} \Biggl\{ n^{-1} \sum_{i=1}^n g^r(X_i; \theta_0)\Biggr\}^2 +
O_p(n^{-3/2}).
\]
Hence,
\begin{eqnarray*}
\Delta_n(\theta_0; a_n) &=& \Delta_n(0) -
2a \sum_{r=1}^{p} \Biggl\{ n^{-1} \sum_{i=1}^n g^r(X_i; \theta_0)\Biggr\}^2\\
&=&
n \{R_1 + R_2 + R_{3a}\}^T \{R_1 + R_2 + R_{3a}\} + O_p(n^{-3/2}),
\end{eqnarray*}
where
\[
R_{3a} = R_3 - \frac{a}{n} R_1.
\]
This proves the first part of Theorem \ref{thm2}.

Again, according to Chen and Cui (\citeyear{CC07}), $R_1 + R_2 + R_3$
have cumulants such that $(1 - b/n) \Delta_n(\theta_0; 0)$ is
approximated by $\chi_p^2$ to $n^{-2}$ precision. Taking advantage of
their proof and using a similar derivation to the proof of Theorem
\ref{thm1}, we find $\Delta_n(\theta_0; a_n)$ with $a_n = b/2 + O_p(
n^{-1/2} )$ is approximated by $\chi_p^2$ to $n^{-2}$ precision. This
completes the proof.
\end{pf*}
\end{appendix}

\printaddresses


\begin{thebibliography}{99}

\bibitem[\protect\citeauthoryear{}{1978}]{BG78}
\textsc{Bhattacharya, R. N.} and \textsc{Ghosh, J. K.} (1978).
On the validity of the Edgeworth expansion.
\textit{Ann. Statist.} \textbf{6} 431--451.
\MR{0471142}

\bibitem[\protect\citeauthoryear{}{2002}]{BN02}
\textsc{Brown, B. W.} and \textsc{Newey, W. K.} (2002).
Generalized method of moments, efficient bootstrapping,
and improved inference.
\textit{J. Bus. Econom. Statist.}
\textbf{20} 507--517.
\MR{1945606}

\bibitem[\protect\citeauthoryear{}{1996}]{BE96}
\textsc{Burnside, C.} and \textsc{Eichenbaum, M.} (1996).
Small-sample properties of GMM-based Wald tests.
\textit{J. Bus. Econom. Statist.}
\textbf{14} 294--308.

\bibitem[\protect\citeauthoryear{}{2008}]{CVA08}
\textsc{Chen, J., Variyath, A. M.} and \textsc{Abraham, B.} (2008).
Adjusted empirical likelihood and its properties.
\textit{J. Comput. Graph. Statist.} \textbf{17} 426--443.
\MR{2439967}

\bibitem[\protect\citeauthoryear{}{1993}]{C93}
\textsc{Chen, S. X.} (1993).
On the accuracy of empirical likelihood confidence regions
for linear regression model.
\textit{Ann. Inst. Statist. Math.} \textbf{45} 621--637.
\MR{1252944}

\bibitem[\protect\citeauthoryear{}{2006}]{CC06}
\textsc{Chen, S. X.} and \textsc{Cui, H. J.} (2006).
On Bartlett correction of empirical likelihood
in the presence of nuisance parameters.
\textit{Biometrika} \textbf{93} 215--220.
\MR{2277752}

\bibitem[\protect\citeauthoryear{}{2007}]{CC07}
\textsc{Chen, S. X.} and \textsc{Cui, H. J.} (2007).
On the second-order properties
of empirical likelihood with moment restrictions.
\textit{J. Econometrics} \textbf{141} 492--516.
\MR{2413478}



\bibitem[\protect\citeauthoryear{}{1998}]{C98}
\textsc{Corcoran, S. A.} (1998). Bartlett adjustment of empirical
discrepancy statistics. \textit{Biometrika} \textbf{85} 967--972.

\bibitem[\protect\citeauthoryear{}{1995}]{CDS95}
\textsc{Corcoran, S. A., Davison, A. C.} and
\textsc{Spady, R. H.} (1995).
Reliable inference from empirical likelihood.
Economics Working Paper 10,
Nuffield College, Univ. Oxford.

\bibitem[\protect\citeauthoryear{}{1991}]{DHR91}
\textsc{DiCiccio, T. J., Hall, P.} and \textsc{Romano, J. P.} (1991).
Empirical likelihood is Bartlett-correctable.
\textit{Ann. Statist.} \textbf{19} 1053--1061.
\MR{1105861}

\bibitem[\protect\citeauthoryear{}{2009}]{EO09}
\textsc{Emerson, S.} and \textsc{Owen, A.} (2009).
Calibration of the empirical likelihood method
for a vector mean. \textit{Electron. J. Statist.} \textbf{3}
1161--1192.

\bibitem[\protect\citeauthoryear{}{1990}]{HS90}
\textsc{Hall, P.} and \textsc{La Scala, B.} (1990).
Methodology and algorithms of empirical likelihood.
\textit{Internat. Statist. Rev.} \textbf{58} 109--127.

\bibitem[\protect\citeauthoryear{}{1996}]{HH96}
\textsc{Hall, P.} and \textsc{Horowitz, J. L.} (1996).
Bootstrap critical values for tests based on
generalized-method-of-moments estimators.
\textit{Econometrica} \textbf{64} 891--916.
\MR{1399222}

\bibitem[\protect\citeauthoryear{}{1982}]{H82}
\textsc{Hansen, L. P.} (1982).
Large sample properties of generalized method of moments estimators.
\textit{Econometrica} \textbf{50} 1029--1054.
\MR{0666123}

\bibitem[\protect\citeauthoryear{}{1997}]{I97}
\textsc{Imbens, G. W.} (1997).
One-step estimators for over-identified generalized method of moments models.
\textit{Rev. Econom. Stud.} \textbf{64} 359--383.
\MR{1456135}

\bibitem[\protect\citeauthoryear{}{1998}]{ISJ98}
\textsc{Imbens, G. W., Spady, R. H.} and \textsc{Johnson, P.} (1998).
Informative theoretic approaches to inference in moment condition models.
\textit{Econometrica} \textbf{66} 333--357.
\MR{1612246}

\bibitem[\protect\citeauthoryear{}{1997}]{KS97}
\textsc{Kitamura, Y.} and \textsc{Stutzer, M.} (1997).
An information-theoretic alternative to generalized method of moments
estimation. \textit{Econometrica} \textbf{65} 861--874.
\MR{1458431}

\bibitem[\protect\citeauthoryear{}{1986}]{LZ86}
\textsc{Liang, K.-Y.} and \textsc{Zeger, S. L.} (1986).
Longitudinal data analysis using generalized linear models.
\textit{Biometrika} \textbf{73} 13--22.
\MR{0836430}

\bibitem[\protect\citeauthoryear{}{1994}]{NM94}
\textsc{Newey, W. K.} and \textsc{McFadden, D.} (1994).
Large sample estimation and hypothesis testing. In
\textit{Handbook of Econometrics} \textbf{4} 2111--2245.
North-Holland, Amsterdam.
\MR{1315971}

\bibitem[\protect\citeauthoryear{}{2004}]{NS04}
\textsc{Newey, W. K.} and \textsc{Smith, R. J.} (2004).
Higher order properties of GMM and generalized empirical likelihood
estimators. \textit{Econometrica} \textbf{72} 219--255.
\MR{2031017}

\bibitem[\protect\citeauthoryear{}{1988}]{O88}
\textsc{Owen, A. B.} (1988). Empirical likelihood ratio confidence
intervals for a single functional. \textit{Biometrika} \textbf{75} 237--249.
\MR{0946049}

\bibitem[\protect\citeauthoryear{}{2001}]{O01}
\textsc{Owen, A. B.} (2001). \textit{Empirical Likelihood}.
Chapman and Hall/CRC Press, New York.

\bibitem[\protect\citeauthoryear{}{1994}]{QL94}
\textsc{Qin, J.} and \textsc{Lawless, J.} (1994). Empirical likelihood
and general equations. \textit{Ann. Statist.} \textbf{22} 300--325.
\MR{1272085}

\bibitem[\protect\citeauthoryear{}{2007}]{S07}
\textsc{Schennach, S. M.} (2007). Point estimation
with exponentially tilted empirical likelihood.
\textit{Ann. Statist.} \textbf{35} 634--672.
\MR{2336862}

\bibitem[\protect\citeauthoryear{}{1997}]{S97}
\textsc{Smith, R. J.} (1997).
Alternative semi-parametric likelihood approaches to
generalized method of moments estimation.
\textit{Economic Journal} \textbf{107} 503--519.

\bibitem[\protect\citeauthoryear{}{2004}]{T04}
\textsc{Tsao, M.} (2004). Bounds on coverage probabilities of the
empirical likelihood ratio confidence regions.
\textit{Ann. Statist.} \textbf{32} 1215--1221.
\MR{2065203}

\end{thebibliography}
\end{document}